\documentclass[a4paper,11pt]{elsarticle}
\makeatletter
\def\ps@pprintTitle{%
 \let\@oddhead\@empty
 \let\@evenhead\@empty
 \def\@oddfoot{\centerline{\thepage}}%
 \let\@evenfoot\@oddfoot}
\makeatother

\usepackage{amssymb,amsmath,mathtools}
\usepackage{cases,resizegather,xcolor,booktabs,multirow}
\usepackage[left=22mm,right=22mm,bottom=30mm]{geometry}
\usepackage[all]{xy}

\usepackage[colorlinks]{hyperref}
\AtBeginDocument{%
  \hypersetup{citecolor=black!70!green,
              linkcolor=black!10!red}}

\usepackage[mathlines]{lineno}

\usepackage{enumitem}
\usepackage{subfigure}
\usepackage{float}
\usepackage[font=normalsize]{caption}

\newtheorem{thm}{Theorem}
\newtheorem{lem}{Lemma}
\newtheorem{cor}{Corollary}

\newtheorem{rmk}{Remark}

\newdefinition{defi}{Definition}
\newproof{prf}{Proof}

\newcommand*{\QED}{\hfill\ensuremath{\square}}

\newcommand{\order}[1]{\left|#1\right|}

\newcommand{\GF}[1]{\mathbb{F}_{#1}}

\begin{document}

\allowdisplaybreaks[4]

\begin{frontmatter}



\title{A Recursive Construction of Permutation Polynomials over $\mathbb{F}_{q^2}$ with Odd Characteristic from R\'{e}dei Functions}


\author[amss,ucas]{Shihui Fu}
\ead{fushihui@amss.ac.cn}
\author[amss]{Xiutao Feng}
\ead{fengxt@amss.ac.cn}
\author[sklois]{Dongdai Lin}
\ead{ddlin@iie.ac.cn}
\author[carleton]{Qiang Wang}
\ead{wang@math.carleton.ca}

\address[amss]{Key Laboratory of Mathematics Mechanization, Academy of Mathematics and Systems Science, Chinese Academy of Sciences, Beijing, 100190, CHINA}
\address[ucas]{School of Mathematical Sciences, University of Chinese Academy of Sciences, Beijing, 100049, CHINA}
\address[sklois]{State Key Laboratory of Information Security, Institute of Information Engineering, Chinese Academy of Sciences, Beijing, 100093, CHINA}
\address[carleton]{School of Mathematics and Statistics, Carleton University, Ottawa, ON K1S 5B6 CANADA}


\begin{abstract}
In this paper, we construct two classes of permutation polynomials over $\mathbb{F}_{q^2}$ with odd characteristic from rational R\'{e}dei functions. A complete characterization of their compositional inverses is also given. These permutation polynomials can be generated recursively. As a consequence, we can generate recursively permutation polynomials with arbitrary number of terms. More importantly, the conditions of these polynomials being permutations are very easy to characterize. For wide applications in practice, several classes of permutation binomials and trinomials are given. With the help of a computer, we find that the number of permutation polynomials of these types is very large. \par
\end{abstract}


\end{frontmatter}


\section{Introduction}

Let $q = p^k$ be the power of a prime number $p$, $\GF{q}$ be a finite field with $q$ elements, and  $\GF{q}[x]$ be the ring of polynomials over $\GF{q}$. We call $f(x) \in \GF{q}[x]$ a \emph{permutation polynomial} over $\GF{q}$ if its associated polynomial mapping $f:c\mapsto f(c)$ from $\GF{q}$ to itself is a bijection. It is well known that every permutation on $\GF{q}$ can be expressed as a permutation polynomial over $\GF{q}$.

Permutation polynomials over finite fields have been a hot topic of study for many years, and have applications in coding theory \cite{TIT:DingH13,FFA:Laigle07}, cryptography\cite{CRYPTO:LidlM83,EL:SchwenkH98,CACM:RivestSA78}, combinatorial designs \cite{JCT:DingY06}, and other areas of mathematics and engineering. For example, Dickson permutation polynomials of order five, i.e., $D_5(x, a)=x^5+ax^3-a^2x$ over finite fields,  led to a 70-year research breakthrough in combinatorics \cite{JCT:DingY06}, gave a family of perfect nonlinear functions for cryptography \cite{JCT:DingY06}, generated good linear codes \cite{TIT:CarletDY05,TIT:YuanCD06} for data communication and storage, and produced optimal signal sets for CDMA communications \cite{TCOM:DingY07}, to mention only a few applications of these Dickson permutation polynomials. For more background material, information about properties, constructions, and applications of permutation polynomials may be found in \cite[Chapter 7]{FFields:LidlN97} and \cite[Chapter 8]{HandbookFFields:MullenP13}. For a detailed survey of open questions and recent results we refer to see \cite{FFA:Hou15}.

Recently, Akbary, Ghioca and Wang derived the following useful criterion in studying permutation functions on finite sets. It first appeared in \cite{FFA:AkbaryGW11} and further developed in \cite{FFA:YuanD11}, \cite{FFA:YuanD14}, among others.

\begin{lem}[AGW Criterion]\label{lemma:AGWCriterion}
Let $A$, $S$ and $\bar{S}$ be finite sets with $\order{S} = \order{\bar{S}}$, and let $f:A\rightarrow A$, $\bar{f}:S\rightarrow \bar{S}$, $\lambda:A\rightarrow S$, and $\bar{\lambda}:A\rightarrow \bar{S}$  be maps such that $\bar{\lambda} \circ f=\bar{f} \circ\lambda$ (see the following commutative diagram).
\[
\xymatrix{
A \ar[r]^{f}\ar[d]^{\lambda} & A \ar[d]^{\bar{\lambda}}\\
S \ar[r]^{\bar{f}} & \bar{S} }
\]

If both $\lambda$ and $\bar{\lambda}$ are surjective, then the following statements are equivalent:
\begin{enumerate}[label=(\roman*)]
\item $f$ is a bijection from $A$ to $A$ (a permuation over $A$);
\item $\bar{f}$ is a bijection from $S$ to $\bar{S}$ and $f$ is injective on $\lambda^{-1}(s)$ for each $s\in S$.
\end{enumerate}
\end{lem}

For any function $f(x)$ defined over $A$, Lemma \ref{lemma:AGWCriterion} can be viewed piecewisely.  Namely,  let $S=\{s_0, s_1, \ldots, s_{\ell-1}\}$, then we have

\[
f(x) =
\left\{
\begin{array}{ll}
f_0(x),             &  \mbox{if }x \in C_0 := \lambda^{-1}(s_0); \\
\quad \vdots        &  \qquad \vdots \\
f_i (x),            &  \mbox{if }x \in C_i := \lambda^{-1}(s_i) ; \\
\quad \vdots        &  \qquad \vdots \\
f_{\ell -1}(x),     &  \mbox{if }x \in C_{\ell-1} := \lambda^{-1}(s_{\ell-1}),
\end{array}
\right.
\]
is a bijection if and only if $\bar{f}$ is a bijection from $S$ to $\bar{S}$ and each $f_i$ is injective over $C_i$ for $0\leq i \leq \ell-1$.

In particular, if we take $A = \GF{q}\setminus\{0\} = \langle\gamma\rangle$, $\ell s = q-1$, $\zeta = \gamma^s$, and $\lambda=\bar{\lambda} = x^s$, then we obtain the following commutative diagram.
\[
\xymatrix{
\GF{q}\setminus\{0\} \ar[r]^{f}\ar[d]^{x^s}         & \GF{q}\setminus\{0\} \ar[d]^{x^s}  \\
S=\{1,\zeta,\dots,\zeta^{\ell-1}\} \ar[r]^{\bar{f}} & S = \{1,\zeta,\dots,\zeta^{\ell-1}\}}
\]

This special case of AGW criterion was obtained earlier in \cite{DM:NiederreiterW05} and \cite{SSC:Wang07}. Namely, the cyclotomic mapping polynomial
\begin{equation}\label{equ:CyclotomicMapping}
P(x) =
\left\{
\begin{array}{ll}
\hspace{3mm} 0,   &  \mbox{if }x=0; \\
A_0 x^r,          &  \mbox{if }x \in C_0 :=\langle \gamma^\ell \rangle \subseteq \GF{q}\setminus\{0\} = \langle \gamma \rangle; \\
\quad\vdots       &  \qquad \vdots \\
A_i x^r,          &  \mbox{if }x \in C_i := \gamma^i C_0; \\
\quad\vdots       &  \qquad \vdots \\
A_{\ell -1} x^r,  &  \mbox{if }x \in C_{\ell-1} := \gamma^{\ell-1} C_0,
\end{array}
\right.
\end{equation}
is a permutation polynomial over $\GF{q}$ if and only if $\gcd(r,s) =1$ and $\bar{f}$ permutes $S= \{1,\zeta,\dots,\zeta^{\ell-1}\}$.

Rewriting it in terms of polynomials leads to the following useful criterion.

\begin{cor}[Park-Lee \cite{AMS:LeeP97}, Wang \cite{SSC:Wang07}, Zieve \cite{PAMS:Zieve09}]\label{cor:BasicCorollary}
Let $q-1 = \ell s$ for some positive integers $\ell$ and $s$. Then $P(x)=x^rf(x^s)$ is a permutation polynomial over $\GF{q}$ if and only if $\gcd(r, s) =1$ and $x^rf(x)^s$ permutes the set $\mu_{\ell}$ of all distinct $\ell$-th roots of unity.
\end{cor}

In this paper, we shall employ this corollary to derive two classes of permutation polynomials over finite fields with odd characteristic. These permutation polynomials can be generated recursively. Moreover, the conditions of these polynomials being permutations are very easy to characterize. As a consequence, we can construct permutation polynomials with arbitrary number of terms. Finally, several classes of permutation binomials and trinomials are given. A complete characterization of their composiyional inverses is also given. Thus, this provides another two class permutation polynomials over finite fields with odd characteristic.

The rest of this paper is organized as follows. In the next section, we recall some basic knowledge about the R\'{e}dei functions. In Section \ref{sec:PermutationPolymials}, two classes of permutation polynomials are presented. As examples, several classes of permutation binomials and trinomials are also given in this section. The compositional inverse functions of these permutation polynomials are discussed in Section \ref{sec:CompositionalInverse}.

\section{R\'{e}dei Functions}\label{sec:RedeiFunctions}

In this section, we present some basic results on the R\'{e}dei functions that will be needed in the sequel. Throughout this paper $\GF{q}$ is the finite field with odd characteristic $p$ having $q$ elements.

Let $n$ be a positive integer, $\alpha\in\GF{q^2}\setminus\{0\}$, then we define the following polynomials over $\GF{q^2}$.
\begin{equation*}
  \begin{split}
    G_n(x,\alpha) & = \sum_{i=0}^{\lfloor\frac{n}{2}\rfloor}\binom{n}{2i}\alpha^ix^{n-2i},\\
    H_n(x,\alpha) & = \sum_{i=0}^{\lfloor\frac{n}{2}\rfloor}\binom{n}{2i+1}\alpha^ix^{n-2i-1}.
  \end{split}
\end{equation*}
The \emph{R\'{e}dei function} is a rational function over $\GF{q^2}$ defined as $R_n(x,\alpha)=\frac{G_n(x,\alpha)}{H_n(x,\alpha)}$. It is easy to check that
\begin{equation}\label{equ:RelationofGH+}
(x+\sqrt{\alpha})^n=G_n(x,\alpha)+ H_n(x,\alpha)\sqrt{\alpha},
\end{equation}
and
\begin{equation}\label{equ:RelationofGH-}
(x-\sqrt{\alpha})^n=G_n(x,\alpha)- H_n(x,\alpha)\sqrt{\alpha}.
\end{equation}
This allows us to generate the polynomials $G_n(x,\alpha)$ and $H_n(x,\alpha)$ recursively,
\[
\begin{split}
G_n(x,\alpha)+H_n(x,\alpha)\sqrt{\alpha} & =(x+\sqrt{\alpha})^n \\
                                         & = (x+\sqrt{\alpha})^{n-1}(x+\sqrt{\alpha})\\
                                         & = \left(G_{n-1}(x,\alpha)+H_{n-1}(x,\alpha)\sqrt{\alpha}\right)(x+\sqrt{\alpha})\\
                                         & = \left(G_{n-1}(x,\alpha)x+H_{n-1}(x,\alpha)\alpha\right)+\left(G_{n-1}(x,\alpha)+H_{n-1}(x,\alpha)x\right)\sqrt{\alpha}.
\end{split}
\]
Therefore, we have
\begin{equation}\label{equ:RecursiveOfN}
\begin{split}
  G_0(x,\alpha) & = 1, \quad H_0(x,\alpha)=0,\\
  G_n(x,\alpha) & = G_{n-1}(x,\alpha)x+H_{n-1}(x,\alpha)\alpha, \\
  H_n(x,\alpha) & = G_{n-1}(x,\alpha)+H_{n-1}(x,\alpha)x,
\end{split}
\end{equation}
or in the matrix form
\begin{equation}\label{equ:RecursiveOfMatrix}
\begin{bmatrix}
G_n(x,\alpha) \\
H_n(x,\alpha)
\end{bmatrix}
=
\begin{bmatrix}
x & \alpha \\
1 & x
\end{bmatrix}
\begin{bmatrix}
G_{n-1}(x,\alpha) \\
H_{n-1}(x,\alpha)
\end{bmatrix}
=
\begin{bmatrix}
x & \alpha \\
1 & x
\end{bmatrix}^n
\begin{bmatrix}
1 \\
0
\end{bmatrix}.
\end{equation}

There is a simple relationship between R\'{e}dei functions and the Dickson polynomials of first kind $D_n(x,a)$. For a comprehensive survey on Dickson polynomials, we refer the readers to \cite[Chapter 7]{FFields:LidlN97} and \cite{DICKSON:LidlMT93}. The Dickson polynomial $D_n(x,a)$ over $\GF{q}$ of first kind of degree $n$ in the indeterminate $x$ and with parameter $a\in\GF{q}$ is given as
\[
D_n(x,a)=\sum_{i=0}^{\lfloor\frac{n}{2}\rfloor}\frac{n}{n-i}\binom{n-i}{i}(-a)^i x^{n-2i}.
\]
It is well known that for any $u_1\in\GF{q}$, $u_2\in\GF{q}$ and positive integer $n$, by Waring's formula, it holds that
\[
u_1^n+u_2^n=D_n(u_1+u_2,u_1u_2).
\]
Then by \eqref{equ:RelationofGH+} and \eqref{equ:RelationofGH-}, we have
\begin{equation}\label{equ:DicksonOfG}
\begin{split}
G_n(x,\alpha) & = \frac{1}{2}\left[(x+\sqrt{\alpha})^n+(x-\sqrt{\alpha})^n\right] \\
              & = \frac{1}{2}\ D_n\left(x+\sqrt{\alpha}+x-\sqrt{\alpha},(x+\sqrt{\alpha})(x-\sqrt{\alpha})\right) \\
              & = \frac{1}{2}\ D_n\left(2x,x^2-\alpha\right).
\end{split}
\end{equation}
Moreover, when $n$ is odd, we have
\begin{equation}\label{equ:DicksonOfH}
\begin{split}
H_n(x,\alpha) & = \frac{1}{2\sqrt{\alpha}}\left[(x+\sqrt{\alpha})^n-(x-\sqrt{\alpha})^n\right] \\
              & = \frac{1}{2\sqrt{\alpha}}\left[(x+\sqrt{\alpha})^n+(-x+\sqrt{\alpha})^n\right] \\
              & = \frac{1}{2\sqrt{\alpha}}\ D_n\left(x+\sqrt{\alpha}-x+\sqrt{\alpha},(x+\sqrt{\alpha})(-x+\sqrt{\alpha})\right) \\
              & = \frac{1}{2\sqrt{\alpha}}\ D_n\left(2\sqrt{\alpha},\alpha-x^2\right).
\end{split}
\end{equation}
In fact, $H_n(x,\alpha)$ is a reversed Dickson polynomial (see more details in \cite{FFA:HouMSY09}).

The following well-known lemma is needed. Since the proof is short, we include it here for the reader's convenience.
\begin{lem}\label{lem:CoprimeG&H}
  $\gcd(G_n(x,\alpha),H_n(x,\alpha))=1$.
\end{lem}

\begin{prf}
  Suppose that $\gcd(G_n(x,\alpha),H_n(x,\alpha))\neq 1$. Then by equation \eqref{equ:RelationofGH+}, we have $(x+\sqrt{\alpha})\large\mid G_n(x,\alpha)$ and $(x+\sqrt{\alpha})\mid H_n(x,\alpha)$, which is equivalent to $G_n(-\sqrt{\alpha},\alpha)=0$ and $H_n(-\sqrt{\alpha},\alpha)=0$. Thus, $G_n(-\sqrt{\alpha},\alpha)-\sqrt{\alpha}\cdot H_n(-\sqrt{\alpha},\alpha)=0$.

  However,
  \[
  \begin{split}
      & G_n(-\sqrt{\alpha},\alpha)-\sqrt{\alpha}\cdot H_n(-\sqrt{\alpha},\alpha) \\
    = \ & \sum_{i=0}^{\lfloor\frac{n}{2}\rfloor}\binom{n}{2i}\alpha^i\left(-\sqrt{\alpha}\right)^{n-2i} -\sqrt{\alpha}\cdot\sum_{i=0}^{\lfloor\frac{n}{2}\rfloor}\binom{n}{2i+1}\alpha^i\left(-\sqrt{\alpha}\right)^{n-2i-1} \\
    = \ & \sum_{i=0}^{\lfloor\frac{n}{2}\rfloor}\binom{n}{2i}(-\sqrt{\alpha})^{2i}\left(-\sqrt{\alpha}\right)^{n-2i} +(-\sqrt{\alpha})\cdot\sum_{i=0}^{\lfloor\frac{n}{2}\rfloor}\binom{n}{2i+1}(-\sqrt{\alpha})^{2i}\left(-\sqrt{\alpha}\right)^{n-2i-1} \\
    =\ & (-\sqrt{\alpha})^n\left(\sum_{i=0}^{\lfloor\frac{n}{2}\rfloor}\binom{n}{2i}+\sum_{i=0}^{\lfloor\frac{n}{2}\rfloor}\binom{n}{2i+1}\right) \\
    = \ & (-\sqrt{\alpha})^n\cdot 2^n\ne 0,
  \end{split}
  \]
  which is a contradiction. The conclusion then follows. \QED
\end{prf}

\section{Two Classes of Permutation Polynomials over $\GF{q^2}$}\label{sec:PermutationPolymials}

In this section, we give the constructions of two classes of permutation polynomials derived from polynomials $H_n(x,\alpha)$ and $G_n(x,\alpha)$. The characterization of being permutations of these constructions are very simple. Besides, from last section we can see that these permutation polynomials can be generated recursively. This, of course, is an advantage in practical applications. Finally, as corollaries, we consider some special cases and construct some classes of permutation binomials and trinomials.

\subsection{General Constructions}

In this subsection, we give the general constructions of two classes of permutation polynomials over $\GF{q^2}$ with odd characteristic, which are presented in the next two theorems.

\begin{thm}\label{thm:ConstructionOfH}
  Suppose $n>0$ and $m$ are two integers. Let $\alpha\in\GF{q^2}$ satisfy $\alpha^{q+1}=1$. Then polynomial
  \[
  P(x)=x^{n+m(q+1)}H_n(x^{q-1},\alpha)
  \]
  permutes $\GF{q^2}$ if and only if any one of the following conditions holds
  \begin{enumerate}
    \item when $\sqrt{\alpha}\in \mu_{q+1}$, $\gcd(n(n+2m),q-1)=1$.
    \item when $\sqrt{\alpha}\notin \mu_{q+1}$, $\gcd(n+2m,q-1)=1$ and $\gcd(n,q+1)=1$.
  \end{enumerate}
\end{thm}

\begin{prf}
  By Corollary \ref{cor:BasicCorollary}, it is equivalent to prove that $\gcd(n+m(q+1),q-1)=1$ and
  \[
  x^{n+m(q+1)}H_n(x,\alpha)^{q-1}
  \]
  permutes $\mu_{q+1}$. Below we assume that $\gcd(n+m(q+1),q-1)=1$, or equivalently $\gcd(n+2m,q-1)=1$. Note that this implies $n$ is odd since $q$ is odd.

  Firstly, we show that $H_n(x,\alpha)$ has no roots in $\mu_{q+1}$. Otherwise, suppose that there exist some $b\in\mu_{q+1}$ such that $H_n(b,\alpha)=0$. Then raising $H_n(b,\alpha)$ to the power of $q$ and substituting $b^q=b^{-1}$ and $\alpha^q=\alpha^{-1}$ into it, we can obtain
  \[
  \begin{split}
  H_n(b,\alpha)^q & = \left(\sum_{i=0}^{\lfloor\frac{n}{2}\rfloor}\binom{n}{2i+1}\alpha^ib^{n-2i-1}\right)^q=\sum_{i=0}^{\lfloor\frac{n}{2}\rfloor}\binom{n}{2i+1}\alpha^{-i}b^{-n+2i+1}. \\
  \end{split}
  \]
  Since $n$ is odd, we can denote $2i'=n-(2i+1)$. The above equation can be reduced to
  \begin{equation}\label{equ:QthPowerOfH}
  \begin{split}
  H_n(b,\alpha)^q & = b^{-n}\alpha^{-\frac{n-1}{2}}\sum_{i'=0}^{\lfloor\frac{n}{2}\rfloor}\binom{n}{2i'}\alpha^{i'}b^{n-2i'}=b^{-n}\alpha^{-\frac{n-1}{2}}G_n(b,\alpha). \\
  \end{split}
  \end{equation}
  This implies that $G_n(b,\alpha)=0$, which contradict to Lemma \ref{lem:CoprimeG&H}. Thus, $H_n(x,\alpha)$ has no roots in $\mu_{q+1}$, and similarly, $G_n(x,\alpha)$ has no roots in $\mu_{q+1}$ as well. Hence, we have $H_n(x,\alpha)\in\GF{q^2}^*$ for any $x\in\mu_{q+1}$. It then implies that $x^{n+m(q+1)}H_n(x,\alpha)^{q-1}\in\mu_{q+1}$ for any $x\in\mu_{q+1}$. In the following, we show that $x^{n+m(q+1)}H_n(x,\alpha)^{q-1}$, or equivalently $x^nH_n(x,\alpha)^{q-1}$ is injective over $\mu_{q+1}$.

  Suppose that there exist $x,y\in\mu_{q+1}$ such that $x^nH_n(x,\alpha)^{q-1}=y^nH_n(y,\alpha)^{q-1}$, then by \eqref{equ:QthPowerOfH}, we can derive that
  \begin{equation}\label{equ:OriginalForm}
  x^nH_n(x,\alpha)^{q-1}=x^n\frac{H_n(x,\alpha)^q}{H_n(x,\alpha)}=\alpha^{-\frac{n-1}{2}}\frac{G_n(x,\alpha)}{H_n(x,\alpha)}.
  \end{equation}
  So $x^nH_n(x,\alpha)^{q-1}=y^nH_n(y,\alpha)^{q-1}$ if and only if
  \begin{equation}\label{equ:GH(x)=GH(y)}
  \frac{G_n(x,\alpha)}{H_n(x,\alpha)}=\frac{G_n(y,\alpha)}{H_n(y,\alpha)}.
  \end{equation}

  Adding $\sqrt{\alpha}$ to the both sides of equation \eqref{equ:GH(x)=GH(y)}, then by \eqref{equ:RelationofGH+}, we can get
  \begin{equation}\label{equ:XYH+}
  \frac{(x+\sqrt{\alpha})^n}{H_n(x,\alpha)}=\frac{(y+\sqrt{\alpha})^n}{H_n(y,\alpha)}.
  \end{equation}
  Similarly, Adding $-\sqrt{\alpha}$ to the both sides of the equation \eqref{equ:GH(x)=GH(y)}, then by \eqref{equ:RelationofGH-}, we can get
  \begin{equation}\label{equ:XYH-}
  \frac{(x-\sqrt{\alpha})^n}{H_n(x,\alpha)}=\frac{(y-\sqrt{\alpha})^n}{H_n(y,\alpha)}.
  \end{equation}
  The remaining proof is divided into the following two cases depending on whether $\sqrt{\alpha}\in\mu_{q+1}$ or not.

  \begin{enumerate}[listparindent=\parindent]
    \item $\sqrt{\alpha}\in \mu_{q+1}$. Suppose $x=\sqrt{\alpha}$, then by \eqref{equ:XYH-}, we must have $y=\sqrt{\alpha}=x$. The case of $x=-\sqrt{\alpha}$ is identical. In the following, we always assume that $x\neq\pm\sqrt{\alpha}$ and $y\neq\pm\sqrt{\alpha}$. By equations \eqref{equ:XYH+} and \eqref{equ:XYH-},
        \[
          \left(\frac{x+\sqrt{\alpha}}{y+\sqrt{\alpha}}\right)^n=\frac{H_n(x,\alpha)}{H_n(y,\alpha)}=\left(\frac{x-\sqrt{\alpha}}{y-\sqrt{\alpha}}\right)^n,
        \]
        or equivalently,
        \begin{equation}\label{equ:NthPowerX=Y}
        \left(\frac{x+\sqrt{\alpha}}{x-\sqrt{\alpha}}\right)^n=\left(\frac{y+\sqrt{\alpha}}{y-\sqrt{\alpha}}\right)^n.
        \end{equation}

        Since $\sqrt{\alpha}\in\mu_{q+1}$, we have $(\sqrt{\alpha})^{q}=\frac{1}{\sqrt{\alpha}}$. Then
        \[
        \left(\frac{x+\sqrt{\alpha}}{x-\sqrt{\alpha}}\right)^{q-1}=\left(\frac{x+\sqrt{\alpha}}{x-\sqrt{\alpha}}\right)^q\left(\frac{x-\sqrt{\alpha}}{x+\sqrt{\alpha}}\right) = \frac{\frac{1}{x}+\frac{1}{\sqrt{\alpha}}}{\frac{1}{x}-\frac{1}{\sqrt{\alpha}}}\left(\frac{x-\sqrt{\alpha}}{x+\sqrt{\alpha}}\right)=-1.
        \]
        This is to say,
        \begin{equation}\label{equ:QthPowerX=Y}
          \left(\frac{x+\sqrt{\alpha}}{x-\sqrt{\alpha}}\right)^{q-1}=-1=\left(\frac{y+\sqrt{\alpha}}{y-\sqrt{\alpha}}\right)^{q-1}.
        \end{equation}

         Finally, it is easy to see that $x^nH_n(x,\alpha)^{q-1}$ is injective over $\mu_{q+1}$ if and only if $x=y$, or equivalently $\frac{x+\sqrt{\alpha}}{x-\sqrt{\alpha}}=\frac{y+\sqrt{\alpha}}{y-\sqrt{\alpha}}$. By equations \eqref{equ:NthPowerX=Y} and \eqref{equ:QthPowerX=Y}, this is possible if and only if $\gcd(n,q-1)=1$.

    \item $\sqrt{\alpha}\notin \mu_{q+1}$. Then $x\neq\pm\sqrt{\alpha}$ and $y\neq\pm\sqrt{\alpha}$. We still have the following equation
        \begin{equation}\label{equ:NthPowerX=Y2}
        \left(\frac{x+\sqrt{\alpha}}{x-\sqrt{\alpha}}\right)^n=\left(\frac{y+\sqrt{\alpha}}{y-\sqrt{\alpha}}\right)^n.
        \end{equation}

        Note that $\sqrt{\alpha}\notin \mu_{q+1}$, we have $(\sqrt{\alpha})^{q}=-\frac{1}{\sqrt{\alpha}}$. Then
        \[
        \left(\frac{x+\sqrt{\alpha}}{x-\sqrt{\alpha}}\right)^{q+1}=\left(\frac{x+\sqrt{\alpha}}{x-\sqrt{\alpha}}\right)^q\left(\frac{x+\sqrt{\alpha}}{x-\sqrt{\alpha}}\right) = \frac{\frac{1}{x}-\frac{1}{\sqrt{\alpha}}}{\frac{1}{x}+\frac{1}{\sqrt{\alpha}}}\left(\frac{x+\sqrt{\alpha}}{x-\sqrt{\alpha}}\right)=-1.
        \]
        This is to say,
        \begin{equation}\label{equ:QthPowerX=Y2}
          \left(\frac{x+\sqrt{\alpha}}{x-\sqrt{\alpha}}\right)^{q+1}=-1=\left(\frac{y+\sqrt{\alpha}}{y-\sqrt{\alpha}}\right)^{q+1}.
        \end{equation}

         Similarly as the case of $\sqrt{\alpha}\in \mu_{q+1}$, $x^nH_n(x,\alpha)^{q-1}$ is injective over $\mu_{q+1}$ if and only if $x=y$, or equivalently $\frac{x+\sqrt{\alpha}}{x-\sqrt{\alpha}}=\frac{y+\sqrt{\alpha}}{y-\sqrt{\alpha}}$. By equations \eqref{equ:NthPowerX=Y2} and \eqref{equ:QthPowerX=Y2}, this is possible if and only if $\gcd(n,q+1)=1$.
  \end{enumerate}
  The proof is now completed. \QED
\end{prf}

\begin{rmk}
  It is worth noticing that in \cite{ARXIV:Zieve13}, the author gave two classes of permutation polynomials over $\GF{q^2}$ by exhibiting classes of low-degree rational functions over $\GF{q^2}$. The permutation polynomials they exhibited can be also represented by R\'{e}dei functions. The first class of permutation polynomials is given as follows.
  \[
  x^{n+m(q+1)}\cdot\left((\gamma x^{q-1}-\beta)^n-\gamma(x^{q-1}-\gamma^q\beta)^n\right),
  \]
  where $\beta, \gamma\in\GF{q^2}$ satisfy $\beta\in\mu_{q+1}$ and $\gamma\notin\mu_{q+1}$. In fact, we can easy verify that the construction in Theorem \ref{thm:ConstructionOfH} is different from the construction presented by Zieve in \cite{ARXIV:Zieve13}. To see this, by equations \eqref{equ:RelationofGH+} and \eqref{equ:RelationofGH-}, we can obtain
  \[
  H_n(x,\alpha)=\frac{(x+\sqrt{\alpha})^n-(x-\sqrt{\alpha})^n}{2\sqrt{\alpha}}.
  \]
  Therefore, the permutation polynomial constructed in Theorem \ref{thm:ConstructionOfH} can be presented in the following form
  \[
    x^{n+m(q+1)}\cdot\frac{\left((x^{q-1}+\sqrt{\alpha})^n-(x^{q-1}-\sqrt{\alpha})^n\right)}{2\sqrt{\alpha}}.
  \]
  A direct observation show that the two constructions are distinct. We can also check that the construction in Theorem \ref{thm:ConstructionOfH} is also distinct from the second class of permutation polynomials constructed by Zieve in \cite{ARXIV:Zieve13}.
\end{rmk}

If we replace $H_n(x,\alpha)$ by $G_n(x,\alpha)$ in Theorem \ref{thm:ConstructionOfH}, we can obtain another class of permutation polynomials over $\GF{q^2}$. Since the proof is almost identical to the proof of Theorem \ref{thm:ConstructionOfH}, we present directly the results as follows and omit the proof.

\begin{thm}\label{thm:ConstructionOfG}
  Suppose $n>0$ and $m$ are two integers. Let $\alpha\in\GF{q^2}$ satisfy $\alpha^{q+1}=1$. Then polynomial
  \[
  P(x)=x^{n+m(q+1)}G_n(x^{q-1},\alpha)
  \]
  permutes $\GF{q^2}$ if and only if any one of the following conditions holds
  \begin{enumerate}
    \item when $\sqrt{\alpha}\in \mu_{q+1}$, $\gcd(n(n+2m),q-1)=1$.
    \item when $\sqrt{\alpha}\notin \mu_{q+1}$, $\gcd(n+2m,q-1)=1$ and $\gcd(n,q+1)=1$.
  \end{enumerate}
\end{thm}

\begin{rmk}
  For a given odd prime power $q$ and an integer $m$, the number of $n$ which satisfies the conditions in Theorem \ref{thm:ConstructionOfH} and \ref{thm:ConstructionOfG} is very large. For instance, when $q=3^k$, $2\leq k\leq 10$ and $m=0$, with the help of a computer, we find that nearly half of the integers in the range $[1,q-1]$ satisfy the condition $\gcd(n(n+2m),q-1)=1$.
\end{rmk}

\subsection{Permutation Binomials and Trinomials}

From last subsection, for any odd $n$, we can construct or generate recursively permutation polynomials with $\frac{n+1}{2}$ terms over $\GF{q^2}$ with odd characteristic. But for particular interesting, permutation polynomials with few terms, especially binomials and trinomials, over finite fields are an active research due to their simple algebraic forms, additional extraordinary properties and their wide applications in many areas of science and engineering. Therefore, in this subsection, we consider some special cases. Some permutation binomials and trinomials are then given.

Let $\gamma$ be a primitive element of $\GF{q^2}$, then $\mu_{q+1}=\langle\gamma^{q-1}\rangle$. An element $\alpha\in\mu_{q+1}$ can be represented as $\alpha=\gamma^{l(q-1)}$, $0\leq l\leq q$, and $\sqrt{\alpha}\in\mu_{q+1}$ if and only if $l$ is even. Then by the recursive relations \eqref{equ:RecursiveOfN} or \eqref{equ:RecursiveOfMatrix}, we have
\[
G_3(x,\alpha)=x^3+3\alpha x,
\]
and
\[
H_3(x,\alpha)=3x^2+\alpha.
\]
By Theorem \ref{thm:ConstructionOfH} and \ref{thm:ConstructionOfG}, we obtain the following permutation binomials over $\GF{q^2}$. Here we may suppose that $3\nmid q$, otherwise, the polynomials constructed by $G_3(x,\alpha)$ and $H_3(x,\alpha)$ will become the monomials, which are permutations over $\GF{q^2}$ if and only if the exponent $d$ is coprime with $q^2-1$.

\begin{cor}\label{cor:BinomialsM}
  Let $q$ be an odd prime power with $3\nmid q$, $m$ and $l$ be two integers. Denote
  \[
    P_1(x)=x^{m(q+1)+3q}+3\gamma^{l(q-1)}x^{m(q+1)+q+2}
  \]
  and
  \[
    P_2(x)=3x^{m(q+1)+2q+1}+\gamma^{l(q-1)}x^{m(q+1)+3}.
  \]
  Then
  \begin{enumerate}
    \item when $l$ is even, the polynomials $P_1(x)$ and $P_2(x)$ permute $\GF{q^2}$ if and only if $\gcd(3(2m+3),q-1)=1$.
    \item when $l$ is odd, the polynomials $P_1(x)$ and $P_2(x)$ permute $\GF{q^2}$ if and only if $\gcd(2m+3,q-1)=1$ and $\gcd(3,q+1)=1$.
  \end{enumerate}
\end{cor}

Specializing even further to the values $m=q-3$, $m=1$, or $m=0$ yields the following consequences.

\begin{cor}\label{cor:MEqualQminusThree}
  Let $q$ be an odd prime power with $3\nmid q$, and $l$ be an integer. Denote
  \[
  P_1(x)=x^{q-2}+3\gamma^{l(q-1)}x^{q^2-q-1}
  \]
  and
  \[
  P_2(x)=3x^{q^2-2}+\gamma^{l(q-1)}x^{q^2-2q}.
  \]
  Then
  \begin{enumerate}
    \item when $l$ is even, the polynomials $P_1(x)$ and $P_2(x)$ permute $\GF{q^2}$ if and only if $q\not\equiv 1\pmod 3$.
    \item when $l$ is odd, the polynomials $P_1(x)$ and $P_2(x)$  permute $\GF{q^2}$ if and only if $q\not\equiv -1\pmod 3$.
  \end{enumerate}
\end{cor}

\begin{rmk}
  Indeed, when $m=q-2$, we still have that $\gcd(2(q-2)+3,q-1)=\gcd(2(q-1)+1,q-1)=\gcd(1,q-1)=1$. Therefore, substituting $m=q-2$ into Corollary \ref{cor:BinomialsM}, we obtain polynomials
  \[
  P_1(x)=x^{2q-1}+3\gamma^{l(q-1)}x
  \]
  and
  \[
  P_2(x)=3x^{q}+\gamma^{l(q-1)}x^{q^2-q+1}.
  \]
  The conclusions in Corollary \ref{cor:MEqualQminusThree} are still valid for the above polynomials $P_1$ and $P_2$.
\end{rmk}

\begin{cor}
  Let $q$ be an odd prime power with $3\nmid q$, and $l$ be an integer. Denote
  \[
  P_1(x)=x^{4q+1}+3\gamma^{l(q-1)}x^{2q+3}
  \]
  and
  \[
  P_2(x)=3x^{3q+2}+\gamma^{l(q-1)}x^{q+4}.
  \]
  Then
  \begin{enumerate}
    \item when $l$ is even, the polynomials $P_1(x)$ and $P_2(x)$ permute $\GF{q^2}$ if and only if $q\not\equiv 1\pmod{3}$ and $q\not\equiv 1\pmod{5}$.
    \item when $l$ is odd, the polynomials $P_1(x)$ and $P_2(x)$ permute $\GF{q^2}$ if and only if $q\not\equiv -1\pmod{3}$ and $q\not\equiv 1\pmod{5}$.
  \end{enumerate}
\end{cor}

\begin{cor}
  Let $q$ be an odd prime power with $3\nmid q$, and $l$ be an even integer. Then the polynomials
  \[
  P_1(x)=x^{3q}+3\gamma^{l(q-1)}x^{q+2}
  \]
  and
  \[
  P_2(x)=3x^{2q+1}+\gamma^{l(q-1)}x^3
  \]
  permute $\GF{q^2}$ if and only if $q\not\equiv 1\pmod{3}$.
\end{cor}

Continuing using the recursive relations \eqref{equ:RecursiveOfN} or \eqref{equ:RecursiveOfMatrix}, we have
\[
G_5(x,\alpha)=x^5+10\alpha x^3+5\alpha^2 x,
\]
and
\[
H_5(x,\alpha)=5x^4+10\alpha x^2+\alpha^2.
\]
By Theorem \ref{thm:ConstructionOfH} and \ref{thm:ConstructionOfG}, we obtain the following permutation trinomials over $\GF{q^2}$.

\begin{cor}\label{cor:TrinomialsM}
  Let $q$ be an odd prime power with $5\nmid q$, $m$ and $l$ be two integers. Denote
  \[
    P_1(x)=x^{m(q+1)+5q}+10\gamma^{l(q-1)}x^{m(q+1)+3q+2}+5\gamma^{2l(q-1)} x^{m(q+1)+q+4}
  \]
  and
  \[
    P_2(x)=5x^{m(q+1)+4q+1}+10\gamma^{l(q-1)}x^{m(q+1)+2q+3}+\gamma^{2l(q-1)}x^{m(q+1)+5}.
  \]
  Then
  \begin{enumerate}
    \item when $l$ is even, the polynomials $P_1(x)$ and $P_2(x)$ permute $\GF{q^2}$ if and only if $\gcd(5(2m+5),q-1)=1$.
    \item when $l$ is odd, the polynomials $P_1(x)$ and $P_2(x)$ permute $\GF{q^2}$ if and only if $\gcd(2m+5,q-1)=1$ and $\gcd(5,q+1)=1$.
  \end{enumerate}
\end{cor}

Similarly, we specialize further to the values $m=q-4$, $m=1$, or $m=0$ obtain the following consequences.

\begin{cor}\label{cor:MEqualQminusFour}
  Let $q$ be an odd prime power with $5\nmid q$, and $l$ be an integer. Denote
  \[
  P_1(x)=x^{2q-3}+10\gamma^{l(q-1)}x^{q^2-2}+5\gamma^{2l(q-1)} x^{q^2-2q}
  \]
  and
  \[
  P_2(x)=5x^{q-2}+10\gamma^{l(q-1)}x^{q^2-q-1}+\gamma^{2l(q-1)}x^{q^2-3q+1}.
  \]
  Then
  \begin{enumerate}
    \item when $l$ is even, the polynomials $P_1(x)$ and $P_2(x)$ permute $\GF{q^2}$ if and only if $q\not\equiv 1\pmod 5$.
    \item when $l$ is odd, the polynomials $P_1(x)$ and $P_2(x)$ permute $\GF{q^2}$ if and only if $q\not\equiv 4\pmod 5$.
  \end{enumerate}
\end{cor}

\begin{rmk}
  Similarly, when $m=q-3$, we still have that $\gcd(2(q-3)+5,q-1)=\gcd(2(q-1)+1,q-1)=\gcd(1,q-1)=1$. Therefore, substituting $m=q-3$ into Corollary \ref{cor:TrinomialsM}, we obtain polynomials
  \[
  P_1(x)=x^{3q-2}+10\gamma^{l(q-1)}x^{q}+5\gamma^{2l(q-1)} x^{q^2-q+1}
  \]
  and
  \[
  P_2(x)=5x^{2q-1}+10\gamma^{l(q-1)}x+\gamma^{2l(q-1)}x^{q^2-2q+2}.
  \]
  The conclusions in Corollary \ref{cor:MEqualQminusFour} are still valid for the above polynomials $P_1$ and $P_2$.
\end{rmk}

\begin{cor}
   Let $q$ be an odd prime power with $5\nmid q$, and $l$ be an integer. Denote
  \[
  P_1(x)=x^{6q+1}+10\gamma^{l(q-1)}x^{4q+3}+5\gamma^{2l(q-1)} x^{2q+5}
  \]
  and
  \[
  P_2(x)=5x^{5q+2}+10\gamma^{l(q-1)}x^{3q+4}+\gamma^{2l(q-1)}x^{q+6}.
  \]
  Then
  \begin{enumerate}
    \item when $l$ is even, the polynomials $P_1(x)$ and $P_2(x)$ permute $\GF{q^2}$ if and only if $q\not\equiv 1\pmod 5$ and $q\not\equiv 1\pmod 7$.
    \item when $l$ is odd, the polynomials $P_1(x)$ and $P_2(x)$ permute $\GF{q^2}$ if and only if $q\not\equiv 4\pmod 5$ and $q\not\equiv 1\pmod 7$.
  \end{enumerate}
\end{cor}

\begin{cor}
   Let $q$ be an odd prime power with $5\nmid q$, and $l$ be an integer. Denote
  \[
    P_1(x)=x^{5q}+10\gamma^{l(q-1)}x^{3q+2}+5\gamma^{2l(q-1)} x^{q+4}
  \]
  and
  \[
    P_2(x)=5x^{4q+1}+10\gamma^{l(q-1)}x^{2q+3}+\gamma^{2l(q-1)}x^{5}.
  \]
  Then
  \begin{enumerate}
    \item when $l$ is even, the polynomials $P_1(x)$ and $P_2(x)$ permute $\GF{q^2}$ if and only if $q\not\equiv 1\pmod 5$.
    \item when $l$ is odd, the polynomials $P_1(x)$ and $P_2(x)$ permute $\GF{q^2}$ if and only if $q\not\equiv 1,4\pmod 5$.
  \end{enumerate}
\end{cor}

Similarly as above permutation binomials and trinomilas, in fact, by the relations \eqref{equ:RecursiveOfN} or \eqref{equ:RecursiveOfMatrix}, we can generate recursively permutation polynomials with arbitrary number of terms. For instance, if $p\nmid \binom{n}{2i}$ for $0\leq i\leq \frac{n-1}{2}$, then both the two classes of permutation polynomials $x^{n+m(q+1)}G_n(x^{q-1},\alpha)$ and $x^{n+m(q+1)}H_n(x^{q-1},\alpha)$ consist of $\frac{n+1}{2}$ terms over $\mathbb{F}_{q^2}$.

\section{Compositional Inverse Functions}\label{sec:CompositionalInverse}

In this section, we consider the compositional inverses of the two classes of permutation polynomials constructed in last section.

Now we return to the cyclotomic mapping polynomial defined as \eqref{equ:CyclotomicMapping}, then the class of permutation polynomials $x^{n+m(q+1)}H_n(x^{q-1},\alpha)$ in Theorem \ref{thm:ConstructionOfH} can be rewritten as
\[
P(x) =
\left\{
\begin{array}{ll}
0,                &  \mbox{if }x=0, \\
A_ix^{n+m(q+1)},  &  \mbox{if }x\in C_i=\{\gamma^{k(q+1)+i}:k=0,1,\dots,q-2\}, 0\leq i\leq q,\\
\end{array}
\right.
\]
where $A_i=H_n(\gamma^{(q-1)i},\alpha)$. Let $\gamma$ be a primitive element of $\mathbb{F}_{q^2}$ and $\omega=\gamma^{q-1}$. Then, for any $x\in C_i$, it is easy to check that
\[
\frac{1}{q+1}\sum_{l=0}^{q}\left(\frac{x^{q-1}}{\zeta^j}\right)^l=\frac{1}{q+1}\sum_{l=0}^{q}\zeta^{l(i-j)}=
\begin{cases}
  0, & \mbox{if } i\neq j, \\
  1, & \mbox{if } i=j.
\end{cases}
\]
Then the class of permutation polynomials $x^{n+m(q+1)}H_n(x^{q-1},\alpha)$ in Theorem \ref{thm:ConstructionOfH} can be uniquely represented as
\[
\begin{split}
P(x) & =\frac{1}{q+1}\sum_{i=0}^{q}A_i x^{n+m(q+1)}\left(\sum_{j=0}^{q}\left(\frac{x^{q-1}}{\zeta^i}\right)^j\right) \\
     & =\frac{1}{q+1}\sum_{i,j=0}^{q}A_i\zeta^{-ij}x^{n+m(q+1)+(q-1)j} \pmod{x^{q^2}-x}.
\end{split}
\]
The compositional inverse can then be determined easily by the Theorem 3.3 in \cite{FFA:ZhengYZP16} or Theorem 1 in \cite{FFA:Wang17}.
\begin{equation}\label{fun:InverseOfH}
P^{-1}(x) =\frac{1}{q+1}\sum_{i,j=0}^{q}\zeta^{t_1 i-(n+m(q+1))ij}\left(\frac{x}{A_i}\right)^{r_1'+(q-1)j} \pmod{x^{q^2}-x},
\end{equation}
where $r_1', t_1\in\mathbb{Z}$ satisfy $(n+m(q+1))r_1'+(q-1)t_1=1$.

Similarly, we can obtain the compositional inverse of the class of permutation polynomials in Theorem \ref{thm:ConstructionOfG}.
\begin{equation}\label{fun:InverseOfG}
P^{-1}(x) =\frac{1}{q+1}\sum_{i,j=0}^{q}\zeta^{t_2 i-(n+m(q+1))ij}\left(\frac{x}{B_i}\right)^{r_2'+(q-1)j} \pmod{x^{q^2}-x},
\end{equation}
where $r_2', t_2\in\mathbb{Z}$ satisfy $(n+m(q+1))r_2'+(q-1)t_2=1$, and $B_i=G_n(\gamma^{(q-1)i},\alpha)$.

From \eqref{fun:InverseOfH} and \eqref{fun:InverseOfG}, both the compositional inverses of the two classes of permutation polynomials consist of $q+1$ terms, and each of the coefficients is a sum of $q+1$ elements. When $q$ is very large, the computing of their compositional inverses is very inefficient. In the following we shall give another form of their compositional inverses, which have a relative simpler form. The following lemma is needed. The original statement is due to Li \emph{et al.} \cite{PREP:LiQW17}, here we make a slight modification to apply to our constructions.

\begin{lem}\label{lem:CompositionalInverses}
  Suppose that $P(x)=x^rf(x^{q-1})$ is a permutation polynomial over $\mathbb{F}_{q^2}$, and $\gcd(r,q^2-1)=1$. Let $I(x)$ be the compositional inverse of the polynomial $x^rf(x)^{q-1}$ over $\mu_{q+1}$, and $r'$ be an integer such that $rr'\equiv 1 \pmod{q^2-1}$. Then the compositional inverse of $P(x)$ over $\mathbb{F}_{q^2}$ is equal to
  \[
  P^{-1}(x)=\left(x^{q^2-q+1}f\left(I\left(x^{q-1}\right)\right)^{q-2}\right)^{r'}I\left(x^{q-1}\right).
  \]
\end{lem}

\begin{prf}
  It is sufficient to show that for any $x\in \mathbb{F}_{q^2}$, it holds that $P^{-1}(P(x))=x$.

  When $x=0$, we can easily see that $P^{-1}(0)=0$. Hence, $P^{-1}(P(0))=P^{-1}(0)=0$.

  When $x\neq 0$, then $P(x)\neq 0$. So we can easily obtain $P(x)^{q-1}\in\mu_{q+1}$, and check that
  \[
  I\left(P(x)^{q-1}\right)=I\left(x^{r(q-1)}f\left(x^{q-1}\right)^{q-1}\right)=I\left(\left(x^{q-1}\right)^rf\left(x^{q-1}\right)^{q-1}\right)=x^{q-1}.
  \]
  Therefore,
  \[
  \begin{split}
  P^{-1}(P(x)) & = \left(P(x)^{q^2-q+1}f\left(I\left(P(x)^{q-1}\right)\right)^{q-2}\right)^{r'}I\left(P(x)^{q-1}\right) \\
               & = \left(x^{r(q^2-q+1)}f\left(x^{q-1}\right)^{q^2-q+1}f\left(x^{q-1}\right)^{q-2} \right)^{r'}x^{q-1} \\
               & = x^{q^2}f\left(x^{q-1}\right)^{r'(q^2-1)} \\
               & = xf\left(x^{q-1}\right)^{r'(q^2-1)}.
  \end{split}
  \]
  Note that $x\neq 0$, then $x^{q-1}\in\mu_{q+1}$. Since $P(x)=x^rf(x^{q-1})$ is a permutation polynomial over $\mathbb{F}_{q^2}$, then by Corollary \ref{cor:BasicCorollary}, we must have that $x^rf(x)^{q-1}$ permutes $\mu_{q+1}$. This implies that $f\left(x^{q-1}\right)\neq 0$. We further obtain that $f\left(x^{q-1}\right)^{q^2-1}=1$. Therefore, $P^{-1}(P(x))=x$.

  This completes the proof. \QED
\end{prf}

Next we computer the compositional inverse function of $x^nH_n(x,\alpha)^{q-1}$ over $\mu_{q+1}$. For the sake of convenience, let $r'$ be an integer such that $(n+m(q+1))r'\equiv 1\pmod{q^2-1}$. Then by equations \eqref{equ:OriginalForm}, \eqref{equ:DicksonOfG} and \eqref{equ:DicksonOfH}, we have
\begin{equation}\label{equ:OriginalEquation}
\begin{split}
x^nH_n(x,\alpha)^{q-1} & = \alpha^{-\frac{n-1}{2}}\frac{G_n(x,\alpha)}{H_n(x,\alpha)} = \alpha^{-\frac{n-1}{2}}\frac{\frac{1}{2}\left((x+\sqrt{\alpha})^n+(x-\sqrt{\alpha})^n\right)}{\frac{1}{2\sqrt{\alpha}}\left((x+\sqrt{\alpha})^n-(x-\sqrt{\alpha})^n\right)} \\
& = \alpha^{-\frac{n}{2}+1}\frac{(x+\sqrt{\alpha})^n+(x-\sqrt{\alpha})^n}{(x+\sqrt{\alpha})^n-(x-\sqrt{\alpha})^n}.
\end{split}
\end{equation}
Below we distinguish two cases according to whether $\sqrt{\alpha}\in\mu_{q+1}$ or not.

\begin{enumerate}[listparindent=\parindent]
\item \label{case:NotinQ+1}
$\sqrt{\alpha}\in\mu_{q+1}$. First we need to verify the hypotheses of Lemma \ref{lem:CompositionalInverses}, namely that $\gcd(n+m(q+1),q^2-1)=1$. This is equivalent to $\gcd(n+m(q+1),q-1)=\gcd(n+2m,q-1)=1$ and $\gcd(n+m(q+1),q+1)=\gcd(n,q+1)=1$. Therefore, we additionally suppose that $\gcd(n,q+1)=1$.

If $x=\sqrt{\alpha}$, then by equation \eqref{equ:OriginalEquation}, we have $(\sqrt{\alpha})^nH_n(\sqrt{\alpha},\alpha)^{q-1}=\alpha^{-\frac{n}{2}+1}$. So the preimage of $\alpha^{-\frac{n}{2}+1}$ is equal to $\sqrt{\alpha}$. Similarly, the preimage of $-\alpha^{-\frac{n}{2}+1}$ is equal to $-\sqrt{\alpha}$.

If $x\neq\pm\sqrt{\alpha}$, then we have
\[
x^nH_n(x,\alpha)^{q-1}=\alpha^{-\frac{n}{2}+1}\frac{\left(\frac{x+\sqrt{\alpha}}{x-\sqrt{\alpha}}\right)^n+1}{\left(\frac{x+\sqrt{\alpha}}{x-\sqrt{\alpha}}\right)^n-1}=l_2(x)\circ x^n\circ l_1(x),
\]
where $l_1(x)=\dfrac{x+\sqrt{\alpha}}{x-\sqrt{\alpha}}$, $l_2(x)=\alpha^{-\frac{n}{2}+1}\dfrac{x+1}{x-1}$ are linear fractional transformations. From equation \eqref{equ:QthPowerX=Y}, we know that $l_1(x)\in\mu_{2(q-1)}$. Note that $\gcd(n,2(q-1))=\gcd(n,q-1)=1$, then $x^n$ permutes $\mu_{2(q-1)}$. Hence we can easily obtain the compositional inverse $I_1(x)$ of function $x^nH_n(x,\alpha)^{q-1}$: $\mu_{q+1}\setminus\{\sqrt{\alpha},-\sqrt{\alpha}\}\rightarrow \mu_{q+1}\setminus\{\alpha^{-\frac{n}{2}+1},-\alpha^{-\frac{n}{2}+1}\}$,
\[
I_1(x)=l_1^{-1}(x)\circ x^{n_1}\circ l_2^{-1}(x) =\sqrt{\alpha}\cdot\frac{\left(\alpha^{\frac{n}{2}-1}x+1\right)^{n_1}+\left(\alpha^{\frac{n}{2}-1}x-1\right)^{n_1}}{\left(\alpha^{\frac{n}{2}-1}x+1\right)^{n_1} -\left(\alpha^{\frac{n}{2}-1}x-1\right)^{n_1}},
\]
where $n_1$ is an integer such that $nn_1\equiv 1\pmod{2(q-1)}$.

Then, a straightforward calculation yields that $I_1(\alpha^{-\frac{n}{2}+1})=\sqrt{\alpha}$, $I_1(-\alpha^{-\frac{n}{2}+1})=-\sqrt{\alpha}$. This implies that the compositional inverse of $x^nH_n(x,\alpha)^{q-1}$ over $\mu_{q+1}$ is equal to $I_1(x)$. By Lemma \ref{lem:CompositionalInverses}, the compositional inverse of $P(x)=x^{n+m(q+1)}H_n(x^{q-1},\alpha)$ over $\mathbb{F}_{q^2}$ is equal to
\[
P^{-1}(x)=x^{r'(q^2-q+1)}H_n\left(I_1(x^{q-1}),\alpha\right)^{r'(q-2)}I_1(x^{q-1}).
\]

\item
$\sqrt{\alpha}\notin\mu_{q+1}$. From $\gcd(n+2m,q-1)=\gcd(n+m(q+1),q-1)=1$ and $\gcd(n,q+1)=\gcd(n+m(q+1),q+1)=1$, we can obtain $\gcd(n+m(q+1),q^2-1)=1$. Noticed that  $\pm\sqrt{\alpha}\notin\mu_{q+1}$, then almost identically to case \ref{case:NotinQ+1}, we can obtain the compositional inverse of $x^nH_n(x,\alpha)^{q-1}$ over $\mu_{q+1}$ is equal to
\[
I_2(x)=l_1^{-1}(x)\circ x^{n_2}\circ l_2^{-1}(x) =\sqrt{\alpha}\cdot\frac{\left(\alpha^{\frac{n}{2}-1}x+1\right)^{n_2}+\left(\alpha^{\frac{n}{2}-1}x-1\right)^{n_2}}{\left(\alpha^{\frac{n}{2}-1}x+1\right)^{n_2} -\left(\alpha^{\frac{n}{2}-1}x-1\right)^{n_2}},
\]
where $n_2$ is an integer such that $nn_2\equiv 1\pmod{2(q+1)}$. In this case the compositional inverse of $P(x)=x^{n+m(q+1)}H_n(x^{q-1},\alpha)$ is equal to
\[
P^{-1}(x)=x^{r'(q^2-q+1)}H_n\left(I_2(x^{q-1}),\alpha\right)^{r'(q-2)}I_2(x^{q-1}).
\]
\end{enumerate}

We know that every permutation on $\GF{q^2}$ can be expressed as a permutation polynomial over $\GF{q^2}$. Indeed we can express the compositional inverse permutation $P^{-1}(x)$ as a polynomial over $\GF{q^2}$. Now let's look at the function $I_1(x)$. By \eqref{equ:OriginalEquation}, we have that
\[
\begin{split}
I_1(x) & = \sqrt{\alpha}\cdot\frac{\left(\alpha^{\frac{n}{2}-1}x+1\right)^{n_1}+\left(\alpha^{\frac{n}{2}-1}x-1\right)^{n_1}}{\left(\alpha^{\frac{n}{2}-1}x+1\right)^{n_1} -\left(\alpha^{\frac{n}{2}-1}x-1\right)^{n_1}} \\
& = \sqrt{\alpha}\cdot \frac{\left(\alpha^{\frac{n-1}{2}}x+\sqrt{\alpha}\right)^{n_1}+\left(\alpha^{\frac{n-1}{2}}x-\sqrt{\alpha}\right)^{n_1}}{\left(\alpha^{\frac{n-1}{2}}x+\sqrt{\alpha}\right)^{n_1} -\left(\alpha^{\frac{n-1}{2}}x-\sqrt{\alpha}\right)^{n_1}} \\
& = \alpha^{\frac{1}{2}}\alpha^{\frac{n_1}{2}-1}\left(\alpha^{\frac{n-1}{2}}x\right)^{n_1}H_{n_1}\left(\alpha^{\frac{n-1}{2}}x,\alpha\right)^{q-1} \\
& = \alpha^{\frac{nn_1-1}{2}}x^{n_1}H_{n_1}\left(\alpha^{\frac{n-1}{2}}x,\alpha\right)^{q-1}.
\end{split}
\]
For the function $I_2(x)$, note that $nn_2\equiv 1\pmod{2(q+1)}$, then we have $\alpha^{\frac{nn_2-1}{2}}=1$. Therefore,
\[
I_2(x)=x^{n_2}H_{n_2}\left(\alpha^{\frac{n-1}{2}}x,\alpha\right)^{q-1}.
\]

We summarize the preceding discussion as the following theorems.

\begin{thm}
  Suppose $n>0$ and $m$ are two integers. Let $\alpha\in\GF{q^2}$ satisfy $\alpha^{q+1}=1$, and $r'$ be an integer such that $\left(n+m(q+1)\right)r'\equiv 1\pmod{q^2-1}$. Define polynomial
  \[
  P(x)=x^{n+m(q+1)}H_n(x^{q-1},\alpha).
  \]
  Then the following statements hold.
  \begin{enumerate}
    \item
    When $\sqrt{\alpha}\in \mu_{q+1}$, $\gcd(n(n+2m),q-1)=1$ and $\gcd(n,q+1)=1$, the compositional inverse function of $P(x)$ over $\GF{q^2}$ is equal to
    \[
    P^{-1}(x)=x^{r'(q^2-q+1)}H_n\left(I_1(x^{q-1}),\alpha\right)^{r'(q-2)}I_1(x^{q-1}),
    \]
    where
    \[
      I_1(x)=\alpha^{\frac{nn_1-1}{2}}x^{n_1}H_{n_1}\left(\alpha^{\frac{n-1}{2}}x,\alpha\right)^{q-1}=\sqrt{\alpha}\cdot\frac{\left(\alpha^{\frac{n}{2}-1}x+1\right)^{n_1}+ \left(\alpha^{\frac{n}{2}-1}x-1\right)^{n_1}}{\left(\alpha^{\frac{n}{2}-1}x+1\right)^{n_1} -\left(\alpha^{\frac{n}{2}-1}x-1\right)^{n_1}},
    \]
    and $n_1$ is an integer such that $nn_1\equiv 1\pmod{2(q-1)}$.

    \item
    When $\sqrt{\alpha}\notin \mu_{q+1}$, $\gcd(n+2m,q-1)=1$ and $\gcd(n,q+1)=1$, the compositional inverse function of $P(x)$ over $\GF{q^2}$ is equal to
    \[
    P^{-1}(x)=x^{r'(q^2-q+1)}H_n\left(I_2(x^{q-1}),\alpha\right)^{r'(q-2)}I_2(x^{q-1}),
    \]
    where
    \[
      I_2(x)=x^{n_2}H_{n_2}\left(\alpha^{\frac{n-1}{2}}x,\alpha\right)^{q-1}=\sqrt{\alpha}\cdot\frac{\left(\alpha^{\frac{n}{2}-1}x+1\right)^{n_2} + \left(\alpha^{\frac{n}{2}-1}x-1\right)^{n_2}}{\left(\alpha^{\frac{n}{2}-1}x+1\right)^{n_2} -\left(\alpha^{\frac{n}{2}-1}x-1\right)^{n_2}},
    \]
    and $n_2$ is an integer such that $nn_2\equiv 1\pmod{2(q+1)}$.
  \end{enumerate}
\end{thm}

Similarly, we can also obtain the compositional inverses of permutation polynomials constructed in Theorem \ref{thm:ConstructionOfG}.

\begin{thm}
  Suppose $n>0$ and $m$ are two integers. Let $\alpha\in\GF{q^2}$ satisfy $\alpha^{q+1}=1$, and $r'$ be an integer such that $\left(n+m(q+1)\right)r'\equiv 1\pmod{q^2-1}$. Define polynomial
  \[
  P(x)=x^{n+m(q+1)}G_n(x^{q-1},\alpha).
  \]
  Then the following statements hold.
  \begin{enumerate}
    \item
    When $\sqrt{\alpha}\in \mu_{q+1}$, $\gcd(n(n+2m),q-1)=1$ and $\gcd(n,q+1)=1$, the compositional inverse function of $P(x)$ over $\GF{q^2}$ is equal to
    \[
    P^{-1}(x)=x^{r'(q^2-q+1)}G_n\left(I_3(x^{q-1}),\alpha\right)^{r'(q-2)}I_3(x^{q-1}),
    \]
    where
    \[
      I_3(x)=\alpha^{n_1+\frac{nn_1+1}{2}}x^{n_1}G_{n_1}\left(\alpha^{\frac{n+1}{2}}x,\alpha\right)^{q-1} = \sqrt{\alpha}\cdot\frac{\left(\alpha^{\frac{n}{2}}x+1\right)^{n_1}-\left(\alpha^{\frac{n}{2}}x-1\right)^{n_1}}{\left(\alpha^{\frac{n}{2}}x+1\right)^{n_1} +\left(\alpha^{\frac{n}{2}}x-1\right)^{n_1}},
    \]
    and $n_1$ is an integer such that $nn_1\equiv 1\pmod{2(q-1)}$.
    \item
    When $\sqrt{\alpha}\notin \mu_{q+1}$, $\gcd(n+2m,q-1)=1$ and $\gcd(n,q+1)=1$, the compositional inverse function of $P(x)$ over $\GF{q^2}$ is equal to
    \[
    P^{-1}(x)=x^{r'(q^2-q+1)}G_n\left(I_4(x^{q-1}),\alpha\right)^{r'(q-2)}I_4(x^{q-1}),
    \]
    where
    \[
      I_4(x)=\alpha^{n_2+1}x^{n_2}G_{n_2}\left(\alpha^{\frac{n+1}{2}}x,\alpha\right)^{q-1} = \sqrt{\alpha}\cdot\frac{\left(\alpha^{\frac{n}{2}}x+1\right)^{n_2}-\left(\alpha^{\frac{n}{2}}x-1\right)^{n_2}}{\left(\alpha^{\frac{n}{2}}x+1\right)^{n_2} +\left(\alpha^{\frac{n}{2}}x-1\right)^{n_2}},
    \]
    and $n_2$ is an integer such that $nn_2\equiv 1\pmod{2(q+1)}$.
  \end{enumerate}
\end{thm}

\section*{Acknowledgments}

This work was supported by the National Key Research and Development Program of China (No. 2016YFB0800401), the National Natural Science Foundation of China (No. 61572491 and 11688101) and NSERC of Canada.


{\footnotesize
\bibliographystyle{alpha}
\bibliography{E:/Code/TeX/BibFiles/BibFile}}

\end{document}